# THE MATHEMATICS

Javier F. A. Guachalla H.

La Paz – Bolivia
2005

**To my family**

# CONTENTS

**PROLOGUE** [1]

The mathematics seen in the context of human knowledge has followed a development not exempted of success and frustration. Success whenever the science consolidated itself giving answers to problems of its study. Frustrations when the method did not properly function, or even more when contradictions were found. Relationship that in their constant dialectic configured the way of what the XIX c. would come to see, mathematics consolidated as an unified science, the rigor properly understood and the XX c. with the object of study and the method duly understood, the rigor completely accepted and implemented.

In the present time,

- In the school education, language and mathematics have become fundamental subjects of the curricula, probably, because they are conceptual basis of it.
- In general there is no scientific career, which has no mathematics in its program of studies.
- Considering the basic sciences in general, we think that a country, particularly in development, which leaves out the basic sciences, will have to postpone technological transference, growing on its dependence.
- It is well understood that the information developed by human knowledge grows more and more rapidly, situation which makes necessary to count with an educational methodology which could make this sustainable. [G2]

The present essay has as objective, to describe the characteristics of the elements of the mathematical science as an area of knowledge, which we consider to be the object of study and methodology of development; the mathematical knowledge in the present time, describing in a short form the diversity that has acquired and the entourage of it; where we understand by entourage, the philosophy of mathematics, the mathematics education and its application.

This work consists of three parts. In the first one, we develop elements which support the philosophical conception of the science which occupies us. We consider the sensorial perception as a generating element of knowledge, independent from language. Then we develop about knowledge itself, in particular the formalization of the mathematical object. Elements about language and the development of logic are given. And a chapter on the philosophy of mathematics with elements as the object of study, the epistemology of it; ending with a consideration on the philosophical problem of the utility of mathematics.

---

[1] The author is candidate to Doctor in mathematics. Faculty in the Schools of Mathematics, Universidad Mayor de San Andrés, (Emeritus year 2000), and Institute Normal Superior "Simón Bolívar".

On the second part, elements of the mathematical science are described, particularly its areas, according to their development in time, showing the diversity has acquired. The characteristics of the mathematical activity at the present time. And the entourage of mathematics, that is, the philosophy of mathematics, with mostly a complement of the first part mentioning briefly the philosophical schools at the end of XIX c. and beginning of the XX c.; elements of applied mathematics, and mathematics education, also.

Finally in part III, a chapter is included on mathematics and development, a vision of mathematics within the cognitive conception of development.

The methodology followed is given by the table of contents, which the author thinks is the appropriated conceptual sequence to reach the objective of this study. At first, separate language from perception to characterize properties of each, then integrate these, distinguishing in the cognitive object differentiated characteristics.

The objective of this work and its methodology have made necessary to consider elements of the theory of knowledge and conceptualize in the framework of philosophy, treatment which probably has a deviation proper of the formation of the author; however, we have tried to maintain a logical and consistent reasoning as much as possible, thus there are differentiated paragraphs which are short, with the intention of making just an affirmation or premise only, obtaining an argument susceptible to be understood in the framework of the common sense, with the objective of sketching the paths established within this science.

The contribution of this work, if any, would be particularly, in the methodology followed and the attempt to formally explain the object of study of mathematics as an object of knowledge. We do not affirm neither a scientific experimental development of the premises made; particularly of those in areas which are not properly mathematics, nor an originality, we see them in the framework of the reference given by experience; situation which has as a result of not counting with an extensive bibliography; however, we have fixed those to which we refer explicitly and some activities which have served as a reference to the author.

This document is particularly addressed to mathematics teachers, and in general, to those persons who in some way or another have an interest or need to know about the conceptualization of this science.

We acknowledge those who contributed in the elaboration of this essay, in particular the Schools of mathematics of the Universidad Mayor de San Andrés and the Institute Normal Superior Simón Bolívar; academic units which promote the mathematical knowledge, and where elements of this research have been developed.

**INTRODUCTION**

- The areas of knowledge

Knowledge thru history has suffered a compartmentalization due to the detail that has acquired and developed in the study of the objects of nature, making the object of study more specialized.

For a better understanding, let us consider knowledge as the rational structure formed from the answer, man has given to the problematic that his relation with the different aspects of his living has presented to him.

Particularly, let us consider three relationships: The relationship of man with himself, which we simplify as the personal relationship. The second, the relation of man with other men, the social relationship. Finally, the relationship of man with its physical environment.

The questioning born from the personal relationship can be considered as physical, psychic, philosophical and other types. According to these we can establish that knowledge has been structured for example in medicine, psychology, philosophy, etc.

The problematic presented to men by their social relationship has conformed knowledge, in for instance: economics, sociology, politics, communication, etc.

From the relation with the environment we have engineering, natural sciences, that is technology in general.

Considering that the answer given by man to his different questionings has a rational nature, we consider logic as the fundament of the rational thinking, Mathematics as fundament of the scientific thinking argued thru this work, and statistics and probabilities as fundament of the stochastic, or probabilistic thinking.

Note that the boundaries between the different parts of knowledge are of a diffuse nature, since they overlap at their boundaries, building terms as physic-mathematician, applied mathematician, etc. However, we try to refer to the nucleus of the different branches of knowledge, which are the essence of them, differentiating them from the others.

- The object of study

A characteristic that an area of knowledge counts with is that it establishes itself as the discipline which studies a type of objects, determines their properties or a certain type of properties, following a particular methodology. Thus for instance, physics studies objects of the nature, trying to determine physical properties, it will have methods, techniques and characteristics which distinguish it from other areas of knowledge.

According to the elements we classify and order, there will be the constitution of an area of knowledge. The compartmentalization of knowledge is produced whenever an object of study is more and more specialized. So it happens that, according to the detail to which we arrive in these activities, new areas of knowledge will develop. It is said that Greeks considered within mathematics, areas as optics, music, astronomy; areas that with time, arrived to be areas of knowledge by themselves. For instance, it is possible that some university has already a career of genetic engineering, career that forty years ago would not exist. Therefore, and definitively, the areas of knowledge suffer a constant dynamic according to their evolution, particularly when parts of them arrive, say, to a maturity enough to constitute by themselves new areas of knowledge.

- Elements of the evolution of mathematics

Thru this work we will refer particularly to five stages of the development of the mathematical thinking, starting from our first ancestors of about 4 million years ago and the Homo Erectus of about 2.5 millions of years ago. The second the Homo Sapiens between 150 to 70 thousand years ago. The third the ancient civilizations till the Greeks, particularly with the formalization of geometry, number systems and the deductive method. Fourth, the renaissance with the implementation of the scientific method. Finally the period from the middle of the XIX c., till the 1930's, which has represented for mathematics, time of conceptual and methodological refoundation.

The first one, characterized particularly by the acquisition of the bipedal position, position that liberated the hands, giving the possibility to handle instruments, and allowed also the evolution of a flexible thoracic chest, evolution that with the Homo Sapiens establishes the articulation of sounds and the language is developed.

The Greek time interests us particularly by the Aristotelian formalization of the deductive method in logic, with the modus ponens as a form of tautological thinking. During the renaissance the scientific method is implemented ".. to guide the reason and discover the truth in sciences.." [D] by Descartes and Bacon. Method that after three centuries has seen a technological growth during the last part of the XX c.; difficult to have been imagined by those that developed it.

The XIX c., has for mathematics a particular significance, it is the century in which the problematic of calculus is resolved with the development of the mathematical analysis. The non Euclidean geometries are discovered, creating in the mathematics community stupor and conceptual conflicts, since Euclidean geometry had been considered for about 20 centuries as an example of theory by its method and the explanation of the geometric physical world. However, with the implementation of the relativity theory in a hyperbolic space, these new geometries defined a new conceptual path of the mathematical world. Finally, in the last quarter of this century, set theory is developed by Cantor, which includes a formalization of the infinite, and determines in the mathematics community a philosophical, methodological excision with the establishment early in the XX c. of the logicist, formalist, intuitionist and constructivist philosophical views of mathematics.

# Part I

# Philosophy of mathematics

# 1. Chapter I. The perception

In this chapter we establish the first elements of the mathematical knowledge, to illustrate somehow, we ask the reader to try to imagine man in prehistory, even before man counted with language. If the reader has difficulty with this, he can instead think of a baby during his first year, before he starts to speak.

We place ourselves in this situation with the intention to isolate the sensorial perception, to consider it, independent from the linguistic structure, determining this way that it is primary, before this one. Then we note that thru sensorial experience we start to "know".

we consider the perception as a result of the sensorial experience, source of a first knowledge. We learn to recognize objects, sounds, colors, etc.; thru a series of activities particularly abstraction, remembering, intuition, recording of memories, finally we learn, knowledge is established, a cognitive structure starts to be build and developed.

We underline the aspect of perception, in which we perceive in general some elements and not all of them. And the process of withdrawing properties, the abstraction.

## 1.1 The form

In perceiving an object, there exist properties that are withdrawn, as for instance the form and color of the object. Forming what we call an image as a result, the action afterwards can be either to forget about it, or try to remember, memorize it. In this case, the process of knowing the object starts. This knowledge will be better according to the assimilating of a larger number of properties of the object of study, which will in general be the result of a larger number of experiences and the skill to memorize.

## 1.2 The Magnitude

If we have a referent of the object that we perceive, it is possible that we would distinguish its magnitude, as an element of comparison with a known one, that is, the possibility to say that it is larger, equal or smaller. We underline, the magnitude as an expression (rational and cognitive) resulting from comparison with an object of reference.

## 1.3 Causality

Causality is one more aspect that we learn from experience, as a sequence of facts, which we will call phenomenon. As for instance, if we approach a hand to fire, we feel the heat, if we expose ourselves to cold, we cloth ourselves, or we may by experience to deduce that in certain situations, if a day is fairly warm, the next it will rain, etc.

This aspect develops a logic in the sense that, given a situation then a consequence is expected to happen. In a first instance, we can talk of a primary logic, which we call intuition, which without a determined method it structures a knowledge which it can arrive to foresee a happening from a given situation, without questioning about the elements of it, that is, without being able to explain clearly, for instance which is the most relevant aspect in the situation, and why is that.

We call then intuition, a primary form of logic, which from the cognitive experience, concludes affirmations, with reasons more or less understandable, not necessarily certain. Let us note the use of the term "certain" in the sense to be demonstrable, without a doubt.

We underline the following aspects. On one hand the sequence of facts as an order of them. Second, the development of a intuitive relation of causality in the phenomenon.

## 1.4 The Induction

Induction is characterized by the fact that it tries to generalize from the knowledge of a particular case. To induce a result can be a risky task, that is why the need of a method which facilitates the possibility of success, minimizing the error as much as possible in the process. In a certain form it is similar to intuition, but it is possibly the result of a need of a more refined method for analysis of the scientific fact, at the same time as it defines it.

### 1.4.1 The scientific method

In the 1630's Descartes within the rationalist school and Bacon within the empiricism, established the scientific method, as m method ".. to guide the reason and find the truth in the sciences ..".

The method is resumed in four points, which are:

" ...
  a. Do not admit anything which is not absolutely evident.
  b. Divide each problem in as many as convenient of particular simpler problems, to solve it in a better way.
  c. Follow in order your thoughts, going from the simplest to the most complex.

   d. Enumerate completely the data of the problem, and go thorough each element of its solution, to make sure that it has been correctly solved.
   ... " [D]

It is to note that the method established this way has been extremely successful, three centuries after, the XX c. has seen a technological development which would certainly difficult to foresee by those that created it on the XVII c.

The scientific method is a methodology in the inductive way of thinking, which tries to minimize the error. However, we can not have certainty when inducing while we do not experiment. And even then we can not have the certainty to assert that if once it has been positive, the next will be so, unless we have knowledge of each and every one of the aspects of the fact.

As a principle, if all the elements of a phenomenon can be repeated, the result should be the same, this would be a principle of causality, without probability in it. The words that we ought take in account are ".. all the elements ..", which is exactly the aspect that in general, in complex phenomena becomes very problematic and possibly in many of them unknown. For instance, in climatology, it is been said, that in theory a flapping of the wings of a bird in one place of the planet can have consequences in the other side of the world as the formation of a tornado.

1.4.2  The complexity of knowledge

The scientific method, as we pointed out in the previous numeral, has been successful as a method of studying the nature and the cohabitation in it by man. However, we must notice that it has also had as a result, a high diversification of knowledge, proof of which we can find at present in phenomena like the globalization, the internet and in every one of the areas of knowledge; for instance in mathematics it is said that no mathematician can now be an universalist, since the knowledge is so vast that a mathematician can not contribute to all the different areas of it, as it still happened at the end of the XIX c., and beginning of the XX c.

1.5  The continuity in perception

To end up the chapter, let us consider one more cognitive concept and its relation with nature. If we see a film in television, we would say without a doubt that the movement is continuous. However, if we think on the old films, on those long rolls of film, around sixty cm. in diameter; these films and in general every film is a sequence of photographies, which are passed 24 to 32 per second, giving the sensation of continuity in the movement. Therefore, we must admit that continuity perceived in a film is the result of our visual sense. Continuity at a first instance, as a result of perception in the world that surround us is a cognitive concept, resulting from sensorial experience; below we will refer to this concept in the framework of mathematics.

## 2 Chapter II. The konwledge

The knowledge as a result of the sensorial experience, is then the assimilation of the drawn information thru abstraction of properties of the objects and phenomena in nature.

### 2.1 The Concept

Thus as a result of perception, a process of learning initiates whose results will be ideas, images, causes and consequences, which as it goes along, knowing them with more or less detail, they will form concepts, significances, i.e. a knowledge, structured by relations of causality, similarity, comparison, etc. This knowledge by classification, ordering, relation and function fundaments the cognitive structure. Note that without the language this knowledge is developed in the framework of intelligence and rationality of the experiential, and phenomenological fact.

In what it follows of the chapter, we argument on the elements which make the cognitive fact, developing the conceptualization.

### 2.2 Association

To understand the object, we try to recognize it in the cognitive structure that we count with. That is, we "search" similar objects to distinguish it and determine what is it. We associate the drawn properties with a set of images known to us of objects that among their properties count in particular with those withdrawn, if possible all of them.

Let us see the following example, in which we see for the first time a fruit (it has been determined it is a fruit), unknown to us and which is green, round but the double in size of a lemon. Possibly, once some information has been picked up, as for example the form and color, we try to determine what it is. If it is round, with a green peel, even if it is of larger size, and since we already know lemons, we might conclude that it is a lemon, in this case a big one, in particular we may deduce that its flavor is acid. That is, the information we have drawn has been contrasted to the knowledge, associating it to a set of known objects.

Let us note here that an element that can be misleading is to deduce with respect to a property which was not in consideration. In the example, to have induced that the unknown fruit was acid can be faulty. As it happens with Tangarines of green peel, which are rather sweet and some of them look like big lemons.

## 2.3 Differentiation

Another operation we carry out as an element of information in the cognitive process is the one of comparing the object, with those in a set of reference.

Comparing as measuring the difference, if any, or determine there is no difference. The difference we determine with respect to an object of reference with one or more of the properties of comparison, which can be qualitative or quantitative.

If the operation of comparison has given as a result that there is no difference, then we induce that we are dealing with the same object or that these two objects are equal, we insist here, with respect to the properties in consideration. In case that there is a difference the operation of comparison can execute an ordering of these objects with respect to the properties being considered.

## 2.4 Integration

Finally, the new information will be integrated to knowledge, however, let us note that this integration as a result of the operations mentioned in the previous numerals, is structural and dynamical, particularly by the characteristics of classification and ordering.

This operation depends on the property under consideration or on the consideration of the object as a whole. For instance, in the previous example we can integrate the object among the green objects, or among fruits, etc., integration which is determined by a greater or lesser detail in the knowledge and the objective of the analysis.

- Classify and order

As a result of his cognitive need, man has established classification and ordering as a methodology to understand the object. Classify in general terms would be to place an object in a set related by properties under consideration. Note that this fact is subjective because of the withdrawing of properties what is tied to the individual capacity but whose objectivity has transcendence in the cognitive fact. And ordering in this situation will be in general the result of a comparison or differentiation; particularly relevant in the knowledge of the natural phenomenon.

- Sets

The concept of a set. We will say that a set is the gathering of objects that have one or more properties in common which distinguish the set as a well defined object, as well as its elements. This definition is one of those primary terms, in which we would need to define with anteriority what is to be understood by "a gathering".

Therefore, to avoid entering into a vicious circle or infinite chain we say that it is a term understood in the framework of common sense, in the understanding that there will not be error of convention in itself.

## 2.5   The cognitive structure

Result of the processes of abstraction, classification, ordering and causality, we structure knowledge, as a reflection of the structure of nature. Let us call just for clarity, scale the detail with what, that knowledge is established, detail which will be relevant or not, in a given moment according to the need.

Note the dynamics of the cognitive structure in the sense that actively determine the scale of the detail of analysis, possibly given by the need of the moment, when we consider the same object as an element of different sets.

## 2.6   Formalization

Formalization as a need in the cognitive process.

We have considered that the process of comparison takes in account properties to be considered; if now we focus attention to the object of reference in the comparison, we notice that if this was going to be performed repeated times, we would consider the "creation" of an object of reference which would serve in different circumstances, and that it counts with properties like being durable, trustable, do not change for instance with climate, and to be reproducible. We are talking then of a unity on one hand and on the other formalizing comparison, because we need so.

This process evolves into the conceptualization of number, as a reflection of quantity, which at the same time is the result of measuring or comparing with respect to a unity. Note the conventional characteristic of a unity.

When the comparison refers to the form of the object instead, the related aspects to this concept will give rise to the need of ordering and classify the elements of it. That is to say it will be established the need to conceptualize the geometrical object.

These two instances in which we have underlined the term need, are examples which give rise to consider the mathematical object as ontologically necessary for science, concept which we will refer to below.

2.7  Knowledge and reality

We end the chapter establishing that knowledge is a structured abstraction of reality, a product of the need of man to understand it and manage himself in it.

This reflection of reality will be as much exact as the details of reality are considered and understood from it.
This aspect leads to ask ourselves, when can we say that we know something? Possibly, it is in this situation that the concept of scale has a major relevance.

## 3 Chapter III. The language

In the first chapter we have isolated the cognitive experience from language to establish particularities of it and differentiate the aspects brought with the apparition of language.

### 3.1 The word

It is considered that around 150 to 70 thousand years ago the Homo Sapiens has started to use the word. It is possible that the first words have been names of objects, then actions, etc. Establishing with time the rules of grammar, the language. Let us think of a baby, he follows this sequence when he starts to talk, first he learns names of persons, objects, etc.

A concept associated to the phenomenon of the apparition of the speech is the significance (semantics), the relationship between the word and the image or concept that represents. In mathematical terms a mapping has been established, meanings and words related by this mapping..

Then if we consider the word, we must admit that it is a name (which we assimilate to the term of being a symbol (audible) with the objective to uniform the analysis), a symbol that represents an object, action or other. With the evolution of this, the language is developed conventionally, with a grammar and its own rules.

We must point out the conventional aspect of language, in the sense that does not exist at the beginning a rule to call for instance an apple the fruit we know as such, more over other persons call it manzana.

### 3.2 The language

We see then language as a result of the evolving capacity of man to articulate sounds, establishing the possibility of communication, and express audibly and then by writing, his knowledge.

Thus we conclude that language is another structure which we will call the linguistic structure, developed by man to communicate, it has a symbolic, representative nature. This structure which form part of the cognitive structure has a different nature from that one that we have referred to in the first chapter, which is rather conceptual, because of what we will in general refer to it as the cognitive structure and the last one as the linguistic structure that becomes in mathematical terms a mapping of the above.

We establish then the linguistic structure as a support of the cognitive structure; several times we first learn the word and then later with experience the concept is complemented. Other times we learn first the concept and then we describe it with words, for instance up and down (the form).

## 3.3  Deduction

The communication as an aspect that has allowed the development of knowledge in the human interrelation establishes through the comprehension the need to structure the discourse, the thought. Then we differentiate in one part the grammatical rule, and the other the formalization of the discourse. Formalization that determines the need of rigor in itself, constructing in the ambiance of logic, a methodology, for instance, calling propositions to those phrases which can take a value of truth, true or false, defining concatenation rules, a negation, thus developing a more advanced logic.

Thus we arrive to the time of Aristoteles in the ancient Greece, who gave the lines of deductive reasoning, from the Modus Ponens. The modus ponens is a tautology which starting with a proposition, let us call it p, and an implication (p implies q) it concludes q, which is the base of the deductive logic.

In this situation if p is true (let us suppose it part of the knowledge) and a theorem proves q from p, then q is also true, and we can incorporate q to knowledge, methodology that is characteristic of mathematics.

## 3.4  Axiomatics

The deductive method acquires a semantic connotation the moment in which the primary propositions are defined, admitted true, from which the theory is developed. This first or primary propositions are called axioms, they must be so evident that they do not need proof. The definition of these axioms has conventional aspects.

## 3.5  True and real

This presents for asking ourselves for the relationship between the implication of reasoning and its correspondence with reality, the true and the real. It is possible to find a right reasoning, even though the correspondence with reality can be more or less clear. This situation is due to the fact in general to the trust and completeness of the data removed by abstraction from the elements which make it, that is, the premises admitted as true.

Thus the cognitive process, counts now with two ways of developing, the perception and the linguistic explanation, however this last one possibly, only as a reference.

The possibility of understanding thru the word or explain thru it, are two processes which they complement each other, in the present time, it is possible that in some schools the word be the factor of a higher development of knowledge, in others it would be the reflection of it.

We point out, that at the beginning of the XX c. the philosophical school of the analytic linguistic developed, among whose representatives was for instance Wittgenstein, this school established that many of the problems of philosophy would probably be due to the deficiency of formality in the language, which in terms of the mapping mentioned above, one of the possible reasons is that this is not biunivocal, what is the intention of formalization.

# 4   Chapter IV. Philosophy of mathematics

To establish the elements of philosophy of mathematics, we complement the development of the first chapters with a summary of the period between XIX c. and the XX c. Period that has definitively influenced the present evolution of this science.

## 4.1   The XIX and XX centuries

In the 1660's the derivative is discovered, developing the calculus, as a very effective tool particularly in the study of physical phenomena. The XVIII c. sees an extraordinary development of modeling in differential equations; however, a fundamental deficiency is felt; Berkeley in the 1730's criticizes the calculus, as a knowledge without a proper theoretical foundation. MacLaurin starts this and it is only with the works of Gauss, Cauchy, Abel, Weiersstrass, Riemann and Dedekind, already in the XIX c. that the mathematical analysis is established as the theoretical foundation of calculus.

On the other hand, during the first half of this century, the non Euclidean geometries were discovered with works by Lowachevski and Boljai. Having Saccheri and Gauss before, discovered elements of these. Discovering that demands an explanation, since during twenty centuries Euclidean geometry had been considered as a model of mathematical knowledge, and as an explanation of the geometric physical space.

The problem of the infinite reappears with the definition of the infinite sets by Cantor, theory that it is not accepted by part of the mathematics community. The analysis has been able to algebrize the problem of limits, and conveniently fundament the work of calculus. However, Cantor's set theory considers, possibly from work of Bolzano, the infinite sets, giving rise to a questioning of a philosophical nature about the concept. It is being questioned if the infinite is accepted or not as a mathematical object. Some rejected it, other accepted it, the mathematics community is divided, it is during this period that the schools of philosophy of mathematics are established, such as the logicism, formalism, intuitionism and constructivism.

This dynamic gives rise to a rethinking of mathematics. On one hand, the problems dealt with in calculus, have in the process of convergence as a central element, the concept of the infinite. On the other hand, the change of an axiom of the Elements of Euclides, gives rise to the development of geometries which are as true as the known one, however, they are difficult to conceptually accept them, until it is found an application of these in the theory of relativity.

Finally the theory of infinite sets is not accepted by part of the mathematics community.

In this situation, the axiomatic deductive method is identified as the method of the mathematical science, then the search for a theory which would be its foundation begins, being determined the set theory as the one, focusing the infinite as a cognitive problem.

However, a step is taken at the time, the mathematical science is of an axiomatic deductive nature, therefore the axioms of it have to be determined. Finally, with the works of Zermelo early in the XX c. and Fraenkel around 1920, the ZFC axiomatic system is established, ZF for Zermelo Fraenkel, C if the axiom of choice is added or not, it is possibly at this time the most accepted system, along with BNG, Bernays, von Neumann, Gödel system, which defines classes to avoid the Russell paradox.

In the 1930's Gödel proves his theorems of incompletitud, theorems that in short, establish that in any mathematical theory, where the arithmetic can be developed, there will be undecidable propositions, in the sense that there will not be possible to prove them true or false. It is interesting to mention that around the same time, in physics, Heisemberg discovers the uncertainty principle that asserts the impossibility of knowing at the same time in quantum physics the velocity and the position in space of a particle. Theoretical results that in some way show some limitation of knowledge.

## 4.2 Ontology

The object of study of mathematics has evolved from the concept of the geometrical form and number at the beginning to the abstract structure in the XX c.

Already at the end of the XIX c. the axioms of the first abstract algebraic structures are established, such as those of groups, vector spaces, and so on. And by the middle of the XX c. the theory of categories is formalized by Cartan, Eilenberg and MacLane, which at present time it is being considered as an ontological base of the mathematical object within the structuralist school.

Following the methodology that we have taken in the first chapters, let us consider that, once the numbers have been formalized as being the abstraction of the concept of quantity, this conception has evolved afterwards while the twenty century passed by into being considered as the elements of what is called now, the structure of the natural numbers.

A set which has certain properties, well established already at the end of the XIX c. by Dedekind and Peano. Several authors have considered the construction of the natural numbers, some considered sequences of the type {Φ}, { Φ,{ Φ}}, etc. where Φ is the empty set. Then, if we carry out the exercise of forgetting the representation and leave the structure, considering the representation as an instantiation of it, then we place ourselves in the domain of structuralism, in which it is asserted the existence of structures, and the distinct representations are merely instantiations, which will be more or less relevant according to the context of the corresponding study.

To establish generalities of the object of study of mathematics, let us first consider the different areas of this and the objects that they study.

- The form. Geometry

In the first chapter we pointed out that among the elements removed by abstraction to know an object, was the form of the object. The study of the form has given rise to the formalization of the objects of study of geometry, as point, straight line, curve, etc.

The generalization, as one of the most important activities in the mathematical science, in geometry has also had several ways of development. For instance, as well as the circle and the sphere have been studied, also the n-dimensional sphere is formalized. The concept of a surface has been generalized to that of a manifold, which are mathematical objects which locally look like Euclidean spaces.

On the other hand we said that to be able to measure, the geometrical object was necessary, making this object to be ontologically necessary for science.

- The quantity. Arithmetic. Algebra

The quantity has giving rise to the abstract concept of number; by the end of the XIX c. numbers have been considered in the framework of set theory, the set of natural numbers, with properties which define the structure of natural numbers.

Likewise, as was mentioned at the beginning of the section, by the same epoch until the beginning of the XX c. the elements necessary to define other algebraic structures have been abstracted, as those of groups, rings, fields, vector spaces, etc. Making the mathematical knowledge structural, at the same time as it established how algebra will be from then on, extending this method to the whole of mathematics.

- Measure. Integration

Measuring as a need in the human doing, has had also in mathematics an evolution towards the abstract element, at the beginning of the XX c. Lebesgue develops measure theory, generalizing the concept of the Riemann integral. Theory that during the XX c. will see a great development in generalization as in applications in technology, as well as the unitification in functional analysis of divers concepts started already in the XVIII c. generalizing characteristics of the vector spaces and linear transformations.

- The infinite

We have already mentioned the infinite; however, referring to the mathematical object, we can not other than to come back to it. What is the infinite? A number? A concept? In complex analysis we talk about the point at infinity, whose adjunction to the complex plane, forms the Riemann sphere, which in a certain way has approached the unbounded component of the plane.

The Greeks had already considered it, and it was during this time, that was object of the paradoxes of Zenon, then they considered two types of the infinite, the potential and the actual one. The potential as for example, the possibility of extending further than any bound the sequence of natural numbers. The actual one as the concept itself.

In Fact, presently its treatment is axiomatic, as for instance the system ZFC considers the axiom of the existence of an infinite set.

- Separation. Topology

At the beginning of the XX c. the concept of a topological space is axiomatized by Hausdorff, axiomatization that abstracts the characteristics of separation and neighborhoods of points in a topological space, the continuity of functions between these spaces becomes a most relevant element in this category. If we make an analogy with physics and chemistry in the study of nature, we would think of geometry and topology having similar roles in mathematics, respectively.

The topological spaces are objects of study of mathematics which count with a high level of abstraction, as well as its theory, they express a structural particularity as a reflection, for instance, the continuity of certain models of nature. Abstract and structural properties which we have insisted in pointing out along the different areas of mathematics.

- Dynamics. Differenciation

Finally we refer to the derivative, a concept developed during the XVII c., theory that was finally, properly conceptualized within the mathematical analysis in the XIX c. and has showed great utility in the phenomenological world. The dynamical systems are a generalization of this problematic, being since the XIX c. the source of a number of theories in mathematics, as well as the abstract structuration of new concepts as for instance, vector fields on topological groups.

### 4.2.1 The abstract structure

We have tried to give in a short form examples in the different areas of mathematics of the object of study, with the end in mind to determine the characteristics of it.

We had considered that the cognitive structure is build thru abstraction and the mathematical object thru formalization. We determine then a particularity which distinguish in the present time the object of study of the mathematical science, it is abstract and structured; result of the abstraction of the structure in knowledge which is at its time the reflection of the structure of nature. The mathematics of nature.

As a conclusion, we say that the mathematical object is the abstract structure, object of knowledge, ontologically necessary to science.

### 4.2.2 Ideal. Real

From this point of view, the object of study of mathematics is then ideal, however, an structured abstraction of the reflection of knowledge of nature. Thus we find two situations, the formalization of the concept in the mathematical object becomes in certain situations to find an ideal object, like the straight line, which can not be found in nature, though it models it. And the other the nature itself, as an structure to be reflected by the knowledge is the real structure.

## 4.3 Epistemology

We have already point out the process of abstraction. After of having abstracted properties of an object, we said that to understand it we classify and order; this process is very general in all the areas of knowledge.

We said in 2.6 that formalization of the unity, and the geometrical object, by cognitive need give rise to the conceptualization of mathematical objects. With the evolution of this science, we can say that the object of mathematics is established in the formalization of the abstract object thru the structural abstraction of the cognitive fact.

In the subtitles that follow we review in a summarized form cognitive activities, characteristic to the science that we study.

### 4.3.1 Structuring

As we saw in the first chapters, we can conclude that classify and order reflect the need to understand the object, the cognitive process. Activities that conform the structure.

We underline that mathematics has incorporated these elements from its method and characteristics, as objects of its knowledge. The activity of classifying has been incorporated in the mathematical knowledge as the relation of equivalence, and ordering as the relation of order.

Mathematics is developed by abstraction of the structure, where structure is been considered in the wide sense of the word, that is, as we pointed out in previous sections, we consider particularities of the elements, functions, relations, the dynamics, etc.

4.3.2 Generalizing

The generalization, with epistemological base in the inductive thought, however we must separate it from it, since in the inductive thought, we try to obtain what we call a law of the phenomenological behavior, while the generalization rather poses us an object whose properties are part of a thesis to be proved in the logical deductive context, and that it can particularize the known object (being generalized).

The generalization is one of the activities in the doing of mathematics generating mathematical knowledge. For instance, the properties of vector space of the Euclidean plane have been generalized to the n-dimensional space, and once these spaces were known, they were generalized to the concept of modules, etc.

The process will most of the times, be done by the axiomatization of the generalizing object, by abstraction of structural properties of the known objects; fact that it is source of theorems to be proved, showing properties about the general object and being able in general to recognize the known object as a particular case.

4.3.3 Deductive logic

From the formalizing of the Modus ponens by Aristoteles and the evolution during the XIX c. the deductive logic is part of the methodology of mathematics, starting with axioms which define the objects of study.

4.3.4 The mathematical induction

The mathematical induction is the result of formalizing the inductive method within this science. It is a theorem in the mathematical theory, therefore it counts with the aspect of mathematical certainty.

It should be noted that this theorem depends on the good ordering principle of the natural numbers, which says "Every nonempty subset of the natural numbers has a minimal element". We ask the reader to stop for a moment and think over the truth of this phrase. This principle has been proved to be equivalent to the axiom of choice, axiom which is one of those of the ZFC system.

### 4.3.5 The mathematical truth

Being the development of the mathematical knowledge logical deductive, and the Modus ponens a tautology, the mathematical knowledge is certain, epistemologically necessary, from its axioms.

## 4.4 Utility of mathematics

It is a problem in present philosophy to explain the utility of mathematics in the real world. Why a mathematical result has an interpretation in the real world?

From the point of view of the analysis being done, we can consider mathematics as a structure of certain knowledge, part of the cognitive structure. Moreover, the knowledge as a reflection of the natural structure gives rise to think that the mathematical theory has been developed from the formalization of the knowledge which reflects the reality is useful in nature, which we call the mathematics of nature. Knowledge that is prior to empirical knowledge in the sense of being ontologically necessary for science.

# Part II

# The mathematics

# 5  Chapter V. The mathematical science

In this chapter we try to describe the mathematical science, starting from its basic areas from a chronological point of view, and also include a description of the characteristics of its actual development.

## 5.1  The basic areas

To design a conceptual mapping of mathematics we adopt a chronological point of view, starting from the experience of the species interpreted in the first chapter already, we consider Geometry, and Algebra as the areas which develop at first. Imagine a square on whose upward corners we write down these two areas of mathematics because of their early appearance (see figure below), that repeating the exercise done with respect to knowledge, they answer the question of which mathematical aspects has man developed in ancient time? So, they evolve and establish themselves, geometry during the time of the Greeks and algebra during the time of Arabians.

The other two corners of the square will be the analysis and the topology, whose appearance or constitution as central areas of mathematics are posterior. The analysis is established from the creation of calculus in the 1600's by Newton and Leibniz and the works of Gauss, Abel, Cauchy, Weierstrass, Riemann and Dedekind in the 1800's developing what will be called the mathematical analysis. The topology starts by the end of the 1800's with the Analysis situs of Poincaré and its axiomatization is published by Hausdorff in 1914.

Finally, we refer to an aspect which has been postponed, and that was developed since the ancient times which is optimization, in the sense that man has always tried to obtain the best result due in principle to his limited nature either to go to a place by the shortest path, or to chase the best prey, etc., aspects that have at the beginning a subjective characteristic, however with the mathematization of it, has been and it is a dynamical part of present mathematics, particularly since the development of calculus, moment from which the determination of maximums and minimums is technically performed. We sketch a rectangle in the center of the square writing in it optimization, activity related to the relation of order. In fact, it can be considered as a transversal area of mathematics, so we find examples in algebra, the partial ordered sets, the lattices, in general the ordered structures; in geometry the geodesics, in analysis the variational calculus, the principles of maximums, we also have results in the framework of convexity etc.

In the paragraphs that follow we mention theories that have been developed within each of the mentioned areas, counting with elements of a conceptual mapping, and to establish in a certain way the diversification of the development that has followed, showing the complexity that the mathematical knowledge has acquired.

A CONCEPTUAL MAP (Resumen)

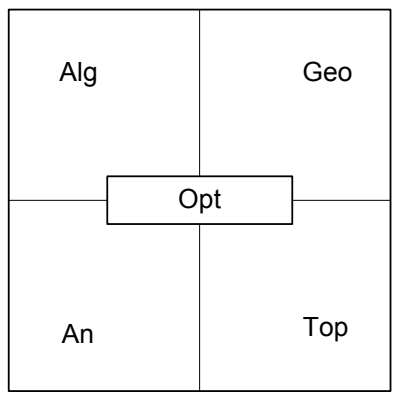

- Geometry

In geometry we point out the Euclidean geometry, the non Euclidean geometries. The differential geometry, the Riemannian geometry. The axiomatic geometry, the algebraic geometry and the differential manifolds, among some theories.

- Algebra

In algebra we can consider the theories related to specific structures like theory of groups, rings, vector spaces, modules, algebras, etc. As other theories like: number theory, Galois theory, Homological algebra and theories that born from the interrelation between areas of mathematics, for instance the geometric algebra. The foundations as logic, set theory and their axiomatic systems, for instance ZFC, and the theory of categories and functors.

- Analysis

In analysis, the real analysis, measure theory, analytic functions or complex analysis, the dynamical systems, control theory. Functional analysis, Operador algebras, etc.

- Topology

In topology, the theories of topological spaces, compact, connected spaces, etc. Degrees of separation. Algebraic topology, differential topology and the theory of topological manifolds. As well as coverings, foliations, fiber bundles, and others.

- Optimization

In optimization, we already named some examples when we saw above the transversal nature of this area in mathematics, we point out some subareas which are relatively new as numerical analysis, operation research of a rather applicative nature.

5.2  Mathematics in the present time and the XXI century

In this subtitle we try in a summarized manner to point out a methodology that mathematics uses more and more in the development of the objects of its study.

In the present time mathematics tends to consider the following elements: Objects, which are generally sets with one or more operations, the structures, i.e. abstract beings with one or more operations or properties. The morphisms, functions with some characteristics which relate structures of one type. This terminology is every time more common inheriting those concepts from theory of categories which summarizes common concepts of different areas of the mathematical knowledge, at the same time that contributes to know itself proving theorems and developing its own theory independently.

We add relations and dynamics, as elements in the study of the mathematical object.

The structures are differentiated by the properties that are subject to study, for instance we have the algebraic structures, the analytical, the topological structures, etc., which are the result of the mathematical work of classification for their study.

The morphisms preserve in general the operations or properties of the type of structure in study, they will depend on the area in which they are considered, for instance the morphisms in topology are the continuous functions, in linear algebra the linear transformations, in the category of groups the homomorphisms of groups, etc.

Among the relations we mention particularly the relations of equivalence, and the relations of order. A relation of equivalence has the particularity of classifying the elements of a set in equivalence classes; it is equivalent to have a partition of a set and to have a relation of equivalence defined on it, giving rise to the quotient set of equivalence classes. These relations are even more important if we consider them in the ambiance of the structures when they preserve the operations making possible to define the quotient structure, thus for instance the quotient group, the quotient vector space etc., are determined.

As far as the relations of order are concerned, they are of a transcendental importance in mathematics as the relations of equivalence. They have the particularity of characterizing the area or activity of optimization which, with the appearance of the computer, take an even more important role in a number of topics, like for instance the operations research, the determination of minimal paths, etc.

Finally, the dynamics as the expression of the evolution of variables with respect to time or other variables, it allows to carry out the study of changes which are obtained from variables with respect to the variation of others, situation whose applicability has been understood already since the creation of calculus.

   The XXI century

The end of the XX c. has seen, may be with some astonishment the development of technology as a result of the consolidation of science and its applications, in possibly most of the areas of knowledge.
As far as mathematics is concerned, it is possible that the topics of non linearity and non commutativity will be the aspects which will develop to some extent during the present century, having the relation of generalization with respect to their predecessors the linearity and commutativity respectively.
The first one already present by the end of the XX c., more over with the advent of the computer, as a tool for obtaining approximations to nature of higher orders than the linear approximation. And the second as a generalization that comes to have the commutative situation as a particular case, presents the possibility to treat a problem with a more unified vision, for which we refer the reader to the work of Alain Connes and his presentation in the International Congress of Mathematics of the 2000 year "The non commutative geometry".

5.3   Paradigms

5.3.1 Discretization and continuity

The advent of the computer has given an impulse to the study of the mathematical problem thru the technique of discretization. Where we understand this as the method of dividing the domain of a problem in a finite number of parts (in general uniformly) obtaining information of a function or system from the information in each of the parts. For instance, in a short form, a technique to calculate the area which is between the graph of a positive function and the x-axis in a given interval, it is possible to approximate the value calculating an average value of the function in each of the subintervals and adding the areas of the rectangles so constructed.

This procedure finds in its study different types of questions like, in how many part it should be divided the domain to obtain a "good" approximation?, that is, with an error less than a predetermined value? Is there such an approximation, or is the process convergent? To what type of functions is the method applicable? etc. method that has given rise to the development of numerical analysis and its application to a every time higher number of areas, in particular to problems of analysis.

On the other hand from the creation of calculus, and its consequences, the mathematical analysis and topology have distinguished the continuity as a fundamental mathematical concept. A situation similar to the theories of waves and corpuscular of light.

5.3.2 Determinism and probability

Going back to the elements developed for the determination of the object of study of mathematics, we have insisted in the abstract nature of it. Thus we say to ourselves that the nature counts with biological, physical, chemical and mathematical characteristics as a cognitive reflection of it.

To complete these facts we can ask ourselves if nature counts with probabilistic characteristics. We will not stop in the philosophical aspects on the determinism of nature. However, we will accept it as a situation of cause and effect which frames the scientific method. And consider that the probability theory is a realization of the scientific method which allows to carry a scientific work counting with a knowledge of the margin of error and the probability, in the certainty of the inexact measure.

We refer to the chaotic characteristic of nature as far as a mathematical expression. This aspect is illustrated by the fact that areas like climatology, or the study of waves in the oceans, are natural phenomena which depend on a very high number of variables and the actual models to which they are assimilated to, result to be very sensitive to the initial conditions, with a chaotic behavior as dynamical systems.

5.4  The entourage of mathematics

In order to complete the analysis, we establish the entourage of the mathematical science as the areas of knowledge with which mathematics has a direct relation. That is, if we go towards the exterior of mathematics, for instance from the activities that a mathematician does, we can consider as the entourage of it, the philosophy of mathematics, the mathematical education, and the application of it, this last one of a very broad nature.

## ENTOURAGE OF MATHEMATICS

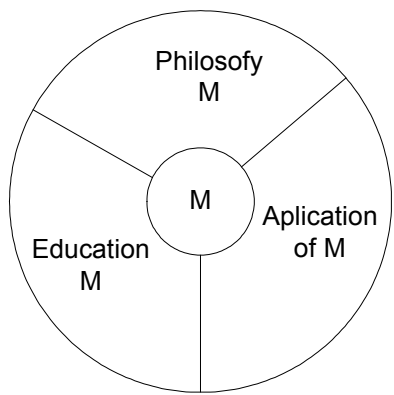

In the following three chapters we try aspects of the entourage of the mathematical science.

# 6 Chapter VI. The philosophical schools

We have referred to the philosophy of mathematics with certain amount of detail in the first part, in this chapter we summarize aspects of the philosophical schools of mathematics established from the evolution of the mathematical knowledge in the XIX c. which would be [BP]:

## 6.1 The logicism

This school affirms that mathematics has the logic as its foundation, from which it is developed, accepts the existence of the infinite set. We can mention among its representatives to Russell.

## 6.2 The formalism

The formalism considers logic as a part of mathematics, that mathematics initiates its study from formal systems, elements given axiomatically, and with logic as an instrument builds the mathematical knowledge, accepts the infinite set. A representative would be Hilbert.

## 6.3 Intuitionism. Constructivism

The intuitionism and the constructivism are schools that accept only the existence of the potential infinite, the constructivists do not accept the principle of the third excluded, therefore their method of demonstration of existence of mathematical objects must be constructive, since the method by contradiction is not valid in this school. Poincaré, Brouwer are representatives of these.

## 6.4 Structuralism

Following the dynamics of the science thru the XX c. we see the appearance of other schools of thinking about the mathematical science, for which we refer the reader to the corresponding bibliography, pointing out particularly the structuralism developed during the second part of the XX c., as the school that fundaments the philosophy of mathematics on the theory of categories.

In the first part of this document, we have tried to develop a structuralist conception of the philosophy of mathematics, trying to explain its evolution from the cognitive need.

# 7 Chapter VII. The applied mathematics

The application of mathematics, is possibly an aspect with a large number of interactions with the development of this science; thru history it has been source of number of mathematical theories, as well as the object of application of others, particularly with the relationship with physics until the XVIII c. At the present time there are more and more areas of knowledge that use mathematical tools, among these we can point out economy, financial mathematics, in biology and social sciences models of population growth, in general those areas of knowledge whose problematic establishes a quantifiable relationship with the variable time, or other variables with dynamical systems, etc., because of this some authors consider mathematics as the language of science.

Mathematics becomes the foundation of the scientific thought by the ontological properties of the mathematical object, this way becomes on one hand a formative instrument in education with the deductive logical thinking, on the other hand as a language or instrument of scientific analysis and conceptualization which is incorporated in the professional doing for the solution of problems in these areas of knowledge, in general expressed as models in terms of mathematical objects. Particularly, since the creation of the derivative and the integral, we count with these instruments which have given rise to a constant growth of dynamical models in the areas of knowledge whose phenomenology is described with elements related to time and/or variables which depend on others, which is possibly a general situation in areas of knowledge structured and quantified.

## 7.1 The mathematical models

Some authors consider the spiral of (empirical) science as: Natural phenomenon / Hypothesis / Solution / Experimentation / Model / Hypothesis / Solution / Experimentation / Model / ….. / Law or natural principle.

As the name suggest, a model is an object that carries information, represents something and its treatment gives back information on the object of study. So we speak of mathematical models which represent a physical or natural phenomenon, in general thru an equation or a system of equations, but not exclusively.

Some examples of areas in which mathematical models are developed in the present time are, for instance agriculture, crops, felling of trees. In the animal kingdom: competition, cohabitation, in natural sciences: physics, chemistry, biology, engineering. In economy: finance mathematics, marketing and others. Some areas of mathematics which are involved in modeling among others are: Ordinary and partial differential equations, integral and functional equations. So much as matrix theory and others.[N]

## 7.2 Nature of the applied mathematics

It is important to point out that the applied mathematics does not differentiate from mathematics, in fact it is important that it counts with all the elements so far established related to the nature of the mathematical knowledge; particularly, structured and logic deductive.

Then if it exists a difference, this can be described saying that while mathematics tries to develop its knowledge by formalization of the abstract structure, the applied mathematics counts with two particular aspects, one the determination of the mathematical model, second the determination of the solution in the context of the mathematical knowledge (see following numeral) to be applied in nature.

Let us note that the applicability of the mathematical concept is not necessarily immediate. For example, around 1860 Weierstrass and Riemann discovered examples of functions that are difficult to imagine intuitively, those which are continuous on their domain and derivable at no point. These have been a century later, contextualized in what is presently called as fractals, which during the last decade have been of application in the technology of compression of information.

## 7.3 The cycle of mathematical modeling

From the point of view of mathematics we consider the cycle of modeling with the following elements.

Phenomenon or problem in nature -> Expression of the problem in the natural language – Expression in a system of equations -> Contextualization in a mathematical theory -> Solution -> Interpretation in the natural language -> Application and/or verification in nature.

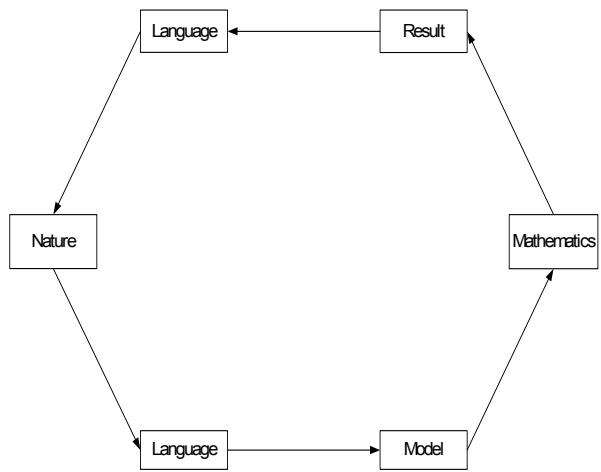

Let us note the language as a vehicle for the comprehension of the problem and the interpretation of the result, as a medium between nature and mathematics.

## 8 Chapter VIII. Mathematical education

For implementing the transmission of knowledge of each of the areas of knowledge is necessary to consider its particularities to learn the principles and methodologies proper to them. The mathematics are not exempted of this fact.

In this chapter we determine a learning cycle for mathematics, besides we consider principles of education as a guide for the analysis. Finally, we refer to the student in the content given by the educational system for the courses of mathematics.

### 8.1 The learning cycle of mathematics

To determine a learning cycle of mathematics we consider as a reference, the cycle of modeling seen above.

The cycle of modeling we have considered is:

Nature -> Language -> Model -> Mathematics -> Solution -> Language -> Nature.

We see that to carry out these elements we must perform the following cycle of activities:

Reflection -> Analysis -> Contextualization -> Deduction -> Interpretation -> Action (praxis).

Which determine the learning cycle of mathematics. Let us called RACLIX (see diagram below).

Reflection: To express in the language or symbolism the problem of nature, we try at this point to understand. It is a divergent activity, creative in the sense to reflect in the cognitive structure to understand, considering if need be different points of view.

Analysis: At this point, once understood, expressed the problem. We analyze it, with the objective in mind to extract the mathematical particularities which formulate the problem in general as a system of equations. At this point, we determine particularly, the data, the unknowns and express the relations among these, in general as equations. It is important to visualize at this point a segregation of the statement into its simplest parts as possible for a better treatment, following the second Descartes principle, and to end up with the mathematic model.

Contextualization: This step possibly immediate from the model, however it calls for a formal treatment of the object, trying to contextualize it in the framework of a mathematical theory, to identify a structure where the elements can count with a number of properties which can be different in another; like the existence of an identity element, inverses, etc., and with the resources of this science, like theorems, for instance, existence and unicity theorems, etc.

Deduction, Logic: As an element of deductive mathematics, which allows from the model to carry out the necessary deductions in the framework of the theory, to arrive to a solution.

Interpretation: Once the solution is determined it will be necessary in general to interpret it in the natural language (or symbolism); it is a form of an inverse operation to analysis, where we have symbolized an object by an unknown; while in this step we interpret the value we found of the unknown as an application to the object.

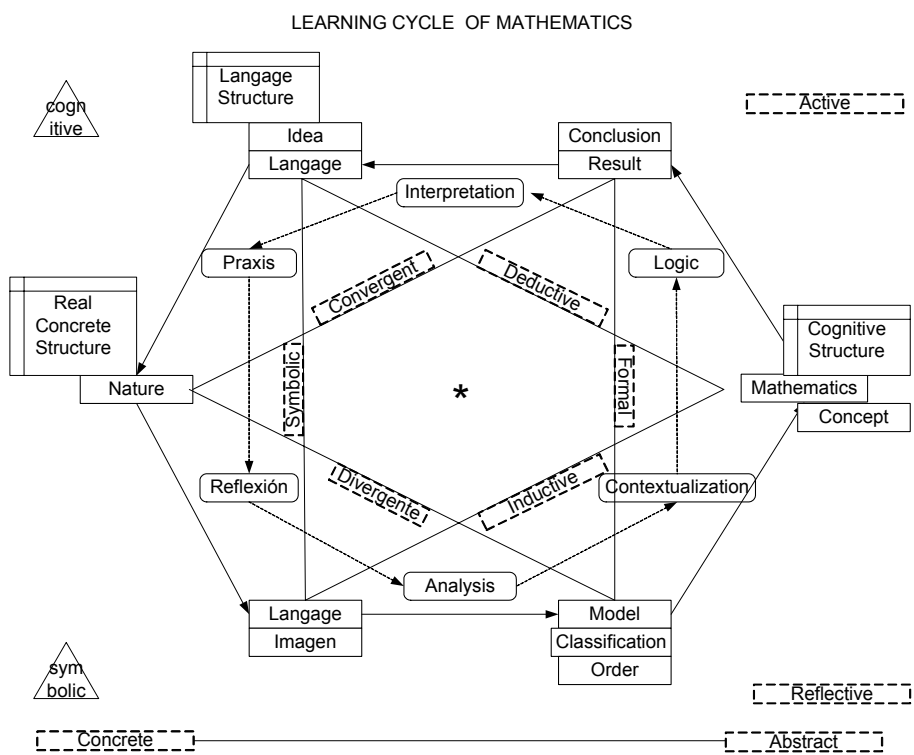

LEARNING CYCLE OF MATHEMATICS

Action, Praxis. We call this way to the doing from the theory. That is, once arrived to this point, we carry out an action to implement the solution as an expression in the nature of the it. This point can also be considered in the framework of a verification. To verify the result in the nature. If the model has not been well determined, it is possible that it is only at this point that an error can be detected even with a very neat solution.

The axis Nature – Mathematics we call it the structural axis, on one hand we have the structure that we have called natural, real, concrete. On the other, the cognitive structure, abstract, as the reflection of the former thru knowledge and in which we emphasize the mathematical structure, resulting from formalism and the deductive method. This axis separates the part above that we characterize as active, while the part below as reflective. We consider the left part as concrete, the right as abstract.

Let us note the language structure (symbolic) backs up the interrelation between the former as an interface. Completing with Mathematics a triangle, that we call symbolic, of symbolic base, representative, expression of the cognitive structure and the mathematics. On the other hand the triangle with base in the model and the solution and vertex in nature, we can consider it a triangle of cognitive base, reflection of the real structure, we call it the cognitive triangle.

On the cognitive triangle we have:

The divergent area: Reflection and analysis. This area tries to reflect the problem into a model, for what it will be necessary to be creative and to follow an analytical methodology.

The formal area: Contextualization and logic. This area is found in the framework of formalism of the mathematical science. Characterized by the deduction and the rigor of the management of the structural mathematical object.

The convergent area: Interpretation and praxis (doing). In this area we take decisions in the framework of the mathematical result, we choose a solution as a function of the context. We converge into doing with the back up of the theory, we carry out the praxis.

In the symbolic triangle, we rather count with:

The symbolic axis: The language as we already explained, as a medium between concrete and abstract.

The inductive area: Analysis and contextualization, represents the experience in the management of abstraction.

The deductive area: Logic and interpretation. Deduction in the framework of the mathematical knowledge and the interpretation as an application of knowledge into reality.

## 8.2 Principles of education

Let us consider the following principles of education:

- Total, partial. Segmented, sequential. Go from the known to the unknown.
- Significative material.
- Transference: Similarity, integration.

Total, partial. Segmented, sequential, from known to unknown

It is to be highlighted that mathematics is essentially sequential, because of the deductive method, which from a premise and an implication it concludes the consequent. Therefore, this aspect is inherent to mathematics.

The curricular grid should be designed to foreseen the advance of courses from the known to the unknown.

However, the segmentation has subjective connotation for the detail of the explanation to which the reasoning arrives, it is an aspect which should be taken in account with certain attention considering the base of knowledge, the experience with the object of study, etc.

The significant learning

The significant learning as Ausbel explains, it establish that for learning to be effective it should be significant, in the sense that it should have a referent for the concept to be learned, to be integrated in the cognitive structure.

Let us note here, that this principle does not enforce the need for the significant to be an element of nature, that is, not necessarily concrete. Particularly, in mathematics, where the abstract object may not have in a first instance other referent than the theory that contextualize it, it is clear however, that a concrete referent is even better, since it enriches itself with the sensorial experience.

In this point paradoxically, the applied mathematics tends to have a certain problem, when because of the developed mathematical knowledge, can not (possibly because of a question of time), detail all the theoretical elements necessary to prove its assertions. In this situation it is convenient to find a coherent explanation, in the framework of the acceptance and consequences of the fact not being proved, definitively it is advisable to give a bibliographic reference, a number of students can be interested on reading about the elements that has not been developed, completing this way their understanding, it is better to name it than to pass it without explanation. And it would be convenient that this happens the lesser number of times as possible during a course.

Transference: Similarity and integration

It is said also that learning should count with elements of similarity and integration with the doing of the subject that learns. This aspect is particularly important in the exemplification and exercising carried out in a course of mathematics, since this science is used for a every time greater number of areas of knowledge, it is of the most importance to pay attention to these two aspects.

8.3  The student in the educational system

The educational system establishes a sequence of studies as primary, secondary and the superior studies. During primary the students take the same courses, in secondary they are divided between social and exact sciences. Finally in the university the students go to their different careers.

It is to be noted that until the university the teaching of mathematics has the character of utilization. Therefore, according to the elements that have been developed in the chapter on applied mathematics, this teaching at the same time it develops the characteristics of the mathematical science, particularly the method logical deductive and the capacity of structural abstraction, it will carry out an exemplification and exercising in an environment rather applicative.

However, in primary school, the pedagogical aspect turns out to be particularly important, since the process of knowing the quantity and the form are initiate, it is at this stage that these elements are established from the natural experience. Therefore we insist in these two aspects analyzed in the first part. On one hand the perception which in principle initiates the elements of form and quantity, independently of language, the symbolism and the logic, the deductive aspect and its backing up to conceptualization.

So we see that there exist three groups of students. One, those who will use mathematics, they are the most, and we referred to this group in the second paragraph of this subtitle. On the other hand we should pay attention to the material of professional study for that group of students who will study to be teachers of mathematics. They need to learn to teach mathematics as a tool (in general) while they should have the necessary formation in mathematics to know its nature and guide the student conveniently in the cognitive process. And the third group, which in number are the less, those who will become mathematicians and work developing mathematics, those who will study mathematics, they will pass from a medium in which the mathematical education is focused in its utilization, to a medium with a focus to develop this science; therefore the first courses will be particularly important with respect to the assimilation of the methodology and objectives of this science.

# Part III

# Mathematics and development

### The strategic character of mathematics in developing countries

# 9 Chapter IX. Mathematics and development

According to the report for the year 2020, the developing countries should invest on the education of their human resources. In this chapter we try to argument that these countries should develop the mathematics education, since we find it strategic for development.

## 9.1 Development

In a first instance we note that the concept of development can have meanings of different nature and deepness. The analysis we carry out in this chapter has the objective to determine the cognitive need in the framework of the present condition.

We consider that we count with two basic components, one the economic the second the educative. We carry out the analysis from the first one merging into the second.

## 9.2 Basic economical activity

We consider as elements of a basic economical activity, particularly the following: the human and the material resources, the infrastructure, the equipment and others. And the matter as input for production, this one as an object to be transformed within production. It is in this last point in which we focus, since to carry out a transformation it is necessary the knowledge.

We generalize the idea to all economical activity, in the sense that in every post of work we have an input and a product.

## 9.3 The cognitive problem

The particularities of knowledge at present should then be established, to consider its treatment in the framework of the decisions to be made.

Thru the whole document we have developed the elements that determine that knowledge has reached in the present time a very high degree of complexity and whose grow is possibly faster and faster. Therefore the developing countries in particular should be able to find a methodology to develop themselves in this environment, where we understand for develop in it: to know it, manage it and produce, produce in the sense of development.

Learn to learn

If the advancement of knowledge is every time faster, the education should find the methodology for the student to count with the possibility to adequate himself. One of the principal characteristic that education has acquired because of this aspect, is that the student must learn to learn.

Where learn to learn would have as objective that the person might in time understand and develop his own knowledge as this one evolves.

Understand and develop it, will be the result of apprehend the new concepts in a structured manner, whose elements detailed in the first part of the document will be an integrated part of the continuous process. It is important then the regular education assists in this labor, counting with bachillers who have the necessary formation for the future.

## 9.4 Mathematics is strategic in developing countries

In this work we have documented the present characteristics of the mathematical science, we mention three, in one part as a science that studies the abstract structure, for other part it borns from the cognitive reflection of the problematic that nature presents us, and its effective utility in the phenomenological world, finally, the methodology logical deductive. Aspects that we resume in the effective utility and the educational characteristics of structuration and logical deductive method.

Learn to structure

If differentiation in knowledge has given rise to complexity, integration complements the process, it is the natural complement which makes viable the cognitive advance. This at its time should ideally be done in a structured manner. In principle the information integrated forms a structure, and it will be the projection of this one in the production which make it transcendent or not.

Mathematics as the science that studies the abstract structure, then takes a relationship of analysis and development with this structuration. The knowledge of this science has as result for one part the treatment of the structure, and for other to count with a effective technological tool, making possible to integrate them by transference in the individual environment.

We conclude that a developing country should consider the mathematical science as an strategic element, to establish a sustainable development in the context of the present knowledge, counting with it as an technological instrument as well as an educational one, in the framework of the characteristics which make mathematics.